\documentclass[11pt]{amsart}
 \usepackage{palatino}
\input{sommers_macros.tex}

\newcommand{\g}{{\mathfrak g}}

\newcommand{\uni}{{\mathfrak u}}

\newcommand{\orbit}{{\mathcal O}}

\newcommand{\al}{{\alpha}}
\newcommand{\coal}{{\alpha \postcheck}}

\newcommand{\complex}{{\mathbf C}}

\newcommand{\zz}{\mathbb Z}

\theoremstyle{definition}

\begin{document}

\setlength{\topmargin}{-.3in}
\setlength{\textheight}{8in} 
\setlength{\oddsidemargin}{0in}
\setlength{\evensidemargin}{-.2in} \setlength{\marginparsep}{0in}
\setlength{\marginparwidth}{0in}

\title {Normality of very even nilpotent varieties in $D_{2l}$}

\author{ Eric Sommers}
\thanks{The author was supported in part by NSF grants DMS-0201826 and DMS-9729992.}
\thanks{The author thanks Viktor Ostrik for directing him to this problem.}

\address{
Department of Mathematics \\
University of Massachusetts--Amherst\\
Amherst, MA 01003}
\address{
Institute for Advanced Study \\
Princeton, NJ 08540}

\date{7/23/03; 9/17/03}

\email{esommers@math.umass.edu}

\begin{abstract}
%%which characteristics:  they also do just complex #'s
For the classical groups, 
Kraft and Procesi \cite{kraft-processi:normal}, \cite{kraft-processi:normal2}
have resolved the question of which nilpotent orbits 
have closures which are normal and which are not,
with the exception of the very even orbits in $D_{2l}$ 
which have partition of the form
$(a^{2k}, b^2)$ for $a, b$ distinct even natural numbers with $a k + b = 2 l$.

In this article, we show that these orbits do have normal closure.  We use the technique
of \cite{sommers:normality}.

\end{abstract}
\maketitle

\section{Some lemmas in $A_l$}

We retain the notation of \cite{sommers:normality}.  
Throughout, $G$ is a connected simple algebraic group over $\complex$, 
$B$ a Borel subgroup, $T$ a maximal torus in $B$.   
The simple roots are denoted by $\Pi$, and they correspond
to the Borel subgroup opposite to $B$.  
Let $\{ \omega_i \}$ be the fundamental weights of $G$ corresponding to $\Pi$.
If $\al \in \Pi$, then 
$P_{\al}$ denotes the parabolic subgroup of semisimple rank one
containing $B$ and corresponding to $\al$.
If $P$ is a parabolic subgroup of $G$, we denote by $\uni_P$ the 
Lie algebra of its unipotent radical.

We recall

\begin{prop} \label{demazure} \cite{demazure:bott-simple}
Let $V$ be a rational representation of $B$ and assume that $V$ extends to 
a representation of the parabolic subgroup $P_{\alpha}$ where $\alpha$ is a simple root.
Let $\lambda \in X^*(T)$ 
be such that $m = \langle \lambda, \coal \rangle \geq -1$.
Then there is a $G$-module isomorphism 
$$ H^i( G/B, V \otimes {\lambda} ) = 
H^{i+1} ( G/B, V \otimes {\lambda \! - \!(m+1)\al} )
\text{ \ for \ all \ } i \in \zz.$$
In particular, if $m=-1$, then all cohomology groups vanish.
\end{prop}

For the rest of this section and the next, let 
$G = SL_{l+1}(\complex)$.  We index the simple roots $\Pi = \{ \al_j \}$
so that $\al_1$ is an extremal root and $\al_j$ is next to $\al_{j+1}$ 
in the Dynkin diagram of type $A_l$.

%Let $P_m$ be the maximal (proper) parabolic subgroup containing $B$	
%corresponding to all the simple roots except $\al_{m}$.  

The following lemma follows easily from 
several applications of the previous proposition.

\begin{lem} \label{repeat:demi1} \cite{sommers:functions}
Let $V$ be a rational representation of $B$ 
which extends to a representation of 
$P_{\al_j}$ for $a \leq j \leq b$.  Let $\lambda \in X^*(T)$ be such that
$\langle \lambda, \al_j \postcheck \rangle = 0$  for $a < j \leq b$.  Set 
$r = \langle \lambda, \al_a \postcheck \rangle$ and assume that
$a-b-1 \leq r \leq -1$.
Then $H^* (V \otimes \lambda) = 0$.  
%The same conclusion holds when 
%$\langle \lambda, \al_j \postcheck \rangle = 0$  for $a \leq j < b$ and
%$r=\langle \lambda, \al_b \postcheck \rangle$.
\end{lem}

A similar statement holds by applying the non-trivial automorphism to the
Dynkin diagram of type $A_{l}$.  
We use this lemma to prove

\begin{lem} \label{repeat:demi2}
Let $V$ be a representation of $B$ 
which is stable under the parabolic subgroups
$P_{\al_j}$ for $1 \leq j \leq b$.  Let $\lambda \in X^*(T)$ be such that
$\langle \lambda, \al_a \postcheck \rangle = 1$ for some $a$ satisfying 
$1 \leq a < b$.
Assume that $\langle \lambda, \al_j \postcheck \rangle = 0$  for $1 \leq j \neq a < b$.
Set $k = \langle \lambda, \al_b \postcheck \rangle$. 
If $-b-1 \leq k \leq -1$ and $k+b-a \neq -1$, then
$H^* (V \otimes \lambda) = 0$.  
\end{lem}

\begin{proof}
If $k+b-a \geq 0$, the result follows directly from Lemma \ref{repeat:demi1}.
On the other hand, if $k+b-a \leq -2$, then as in the proof of Lemma \ref{repeat:demi1}
in \cite{sommers:functions}, 
$$H^i (V \otimes \lambda) = H^{i+b-a}(V \otimes \mu)$$
where 
$$\mu = \lambda + (-k-1)\al_{b} + (-k-2)\al_{b-1} + \dots + (-k-b+a)\al_{a+1}.$$
Now $\langle \mu, \al_j \postcheck \rangle = 0$ for $1 \leq j < a$
and $\langle \mu, \al_a \postcheck \rangle = k+b-a+1$.
By the hypothesis on $\lambda$ and the present assumption about $k+b-a$, we have
$$-a \leq k+b-a+1 \leq -1.$$% (and $-a+1 \leq -1$, meaning $a \neq 1$).
Then Lemma \ref{repeat:demi1} yields the desired vanishing.
\end{proof}

\section{A theorem for $A_l$ (review)}

%%%%notation of max v. min parabolic
Let $P_m$ denote the maximal proper parabolic subgroup of $G = SL_{l+1}(\complex)$
containing $B$ corresponding to all the simple roots except $\al_{m}$.  
Denote the Lie algebra of the unipotent radical of $P_m$ by $\uni_m$.  
The action of $P_m$ on $\uni_m$ gives
a representation of $P_m$ (and also $B$).  Denote the dual representation
by $\uni^*_m$. Set $m' = \mbox{min} \{ m, \  l+1-m \}$.
In \cite{sommers:functions}, Lemma \ref{repeat:demi1} and Proposition \ref{demazure}
were used to prove

\begin{thm} \cite{sommers:functions} \label{shift:A}
Let $r$ be an integer in the range $2m'-2-l \leq r \leq 0$.
Then there is a $G$-module isomorphism 

$$ H^i( G/B, S^n \uni^*_m \otimes r \omega_m ) = 
H^{i} ( G/B, S^{n+rm'} \uni^*_{ \ l+1-m} \otimes -r \omega_{l+1-m} )
\text{ \ for \ all \ } i, n \geq 0.$$
\end{thm}

\section{A theorem for $D_{2l+1}$}

Theorem \ref{shift:A} has an analog in type $D_{2l+1}$. % and in fact, in general...
We label the simple roots of $G$ of type $D_{2l+1}$ as in \cite{onish-vin:book}, 
so $\al_{2l-1}$ lies at the branched vertex of the Dynkin diagram.
Let $P$ be the maximal proper parabolic subgroup 
containing $B$ corresponding to all the simple roots
except $\al_{2l}$.  And let $P'$  
be the maximal proper parabolic subgroup containing $B$
corresponding to all the simple roots except $\al_{2l+1}$ (so $P$ and $P'$
are interchanged by an outer automorphism of $G$). %which fixes T

\begin{thm} \label{shift:D}
Let $r$ be an integer in the range $-3 \leq r \leq 0$.
Then there is a $G$-module isomorphism 

$$ H^i( G/B, S^n \uni^*_P \otimes r \omega_{2l} ) = 
H^{i} ( G/B, S^{n+rl} \uni^*_{P'} \otimes -r \omega_{2l+1} )
\text{ \ for \ all \ } i, n \geq 0.$$
\end{thm}

\begin{proof}

\medskip

{\bf Step 1.}

In this step, $r$ may be an arbitrary integer.
Consider the intersection $V = \uni_P \cap \uni_{P'}$.  
We will show in Step 1  that for all $i, n$
\begin{equation} \label{step1}
H^{i} (S^{n} \uni^*_P \otimes r \omega_{2l} ) = 
H^{i} (S^{n} V^* \otimes r \omega_{2l} ).
\end{equation}

We begin by taking the Koszul resolution of the short exact sequence
$$0 \to U \to \uni^*_P \to V^{*} \to 0$$ (this defines $U$)
and tensoring it with $r \omega_{2l}$.
This gives  
$$ 0 \to \dots \to S^{n-j} \uni^*_P \otimes \wedge^j U \otimes r \omega_{2l} \to \dots 
\to S^n  \uni^*_P \otimes r \omega_{2l} \to S^n V^{*} \otimes r \omega_{2l} \to 0.$$

We claim that $H^{*} (S^{n-j} \uni^*_P \otimes \wedge^{j} U \otimes r \omega_{2l})=0$
for $1 \leq j \leq \dim U$ from which Equation \ref{step1} will follow.
The $T$-weights of $U$ are those 
of the form $\al_k + \al_{k+1} + \dots + \al_{2l}$,
where $1 \leq k \leq 2l$. 
Therefore, if $\lambda$ is a $T$-weight of $\wedge^{j} U$, 
then $\lambda$ is of the form 
$$(0,\dots,0,1,\dots,1,2, \dots, 2, \dots, j-1, \dots, j-1, j, \dots, j, 0)$$ in the basis
of simple roots.  If this expression contains 
a subsequence of the form $m, m, m+1$, 
then $\lambda$ will have inner product $-1$ with the simple coroot
corresponding to the middle $m$.  Hence $H^*(Q \otimes \lambda)= 0$
where $Q$ is any $P$-representation by Proposition \ref{demazure}.
%with $Q:= S^{n-j} \uni^*_P \otimes r \omega_{2l}$.
The same result holds if there are any $0$'s in the initial part of the expression.
Therefore, we are reduced to considering those $\lambda$ of the 
form 
$$(1,2,3, \dots, j-1, j, j, \dots, j, 0).$$
Such a $\lambda$ satisfies 
$\langle \lambda, \al_{2l+1} \postcheck \rangle = \! -j$
with the exception of the case $j=2l$, where instead 
$\langle \lambda, \al_{2l+1} \postcheck \rangle = \! -j+1= \! -2l+1$.
In the latter case $H^{*}(Q \otimes \lambda)=0$ by Lemma \ref{repeat:demi1} applied to 
to the parabolic subgroup with Levi factor of type $A_{2l}$ 
consisting of all simple roots except $\al_{2l}$.
For the cases where $j < 2l$, 
we can apply Lemma \ref{repeat:demi2}, also for the $A_{2l}$ consisting of all simple roots
except $\al_{2l}$.  In that case, $a=j$, $b=2l$, $k=-j$
and so $k + b -a = 2l -2j$, which, being an even number, is never $-1$.  
Also, clearly $-b-1 \leq k \leq -1$.
Thus we conclude that for all weights $\lambda$ appearing in $\wedge^{j} U$,
we have $H^*(Q \otimes \lambda)= 0$ for any $P$-representation $Q$.
Hence for $Q:=S^{n-j} \uni^*_P \otimes r \omega_{2l}$,
it follows that $H^*(Q \otimes \wedge^{j} U)= 0$ by the usual filtration argument.

\medskip

{\bf Step 2.}

Let $V_1$ be the $B$-stable subspace of $\uni$ consisting of the direct sum of
all root spaces $\g_{\al}$
where $-\al$ is bigger than or equal to the root
$$(0, \dots, 0, 1, 2, 1, 1)$$
in the usual partial ordering on roots.
Let $V_2$ be the $B$-stable subspace of $\uni$ consisting of 
the direct sum of all root spaces $\g_{\al}$
where $-\al$ is bigger than or equal to the root
$$(0,0, \dots, 0, 1, 2, 2, 1, 1).$$  
Let $\mu$ be a weight of the form $r \, \omega_{2l} + s \, \omega_{2l+1}$
where $r, s$ are integers.  Assume that $-3 \leq r \leq -1$
and that $s = 0$ if $r=-3$.
In this step we show for all $n \geq 0$ that 
\begin{equation} \label{step2}
H^{*} (S^{n} V_1^* \otimes \mu)=0.
\end{equation}

Take the Koszul resolution of 
$$0 \to U_2 \to V_1^* \to V_2^* \to 0$$
(this defines $U_2$)
and tensor it with $\mu$.  
We will show that 
$$H^*(S^n V_2^* \otimes \mu) = 0$$
and 
$$H^{*} (S^{n-j} V_1^* \otimes \wedge^j U_2 \otimes \mu )=0$$
for  $1 \leq j \leq 2l-2$ and then Equation (\ref{step2}) will follow
(the dimension of $U_2$ is $2l-2$ as shown below).

The subspace $V_2^*$ is stable under 
the minimal parabolic subgroups $P_{\al_m}$ for $m= 2l-1, 2l,$ and $2l+1$.
It follows from the assumption on $\mu$ that
$H^*(S^n V_2^* \otimes \mu) = 0$ by Lemma \ref{repeat:demi1}
applied to the $A_3$ determined by the simple roots
$\al_m$ for $m= 2l-1, 2l,$ and $2l+1$.

Now the $T$-weights of $U_2$ are 
$$\al_k + \al_{k+1} + \dots + \al_{2l-2} + 2 \al_{2l-1} + \al_{2l} +\al_{2l+1}$$ 
where $1 \leq k \leq 2l-2$.
If $\lambda$ is a weight of $\wedge^{j} U_2$, then $\lambda$ is of the form
$$(0,\dots,0,1,\dots,1,2, \dots, j\!-\!1, j, \dots, j, 2j, j, j)$$ in the basis
of simple roots.  As in the previous step, if there are any $0$'s present
or if any of the integers between $1$ and $j-1$ inclusive are repeated,
then $$H^{*}(Q \otimes \lambda)=0$$
where $Q:=S^{n-j} V_1^* \otimes \mu$ since $Q$ is stable 
under the action of the parabolic subgroups $P_{\al_k}$ for $1 \leq k \leq 2l-2$.
Hence we are reduced to considering those $\lambda$ of the form
$$(1, 2, 3, \dots,j\!-\!2, j\!-\!1, j, \dots, j, 2j, j, j)$$
for $1 \leq j \leq 2l-2$.
Such a $\lambda$ satisfies 
$\langle \lambda, \al_{2l-2} \postcheck \rangle = -j$ 
with the exception of $j=2l-2$ where
$\langle \lambda, \al_{2l-2} \postcheck \rangle = -2l+3$.
In the latter case $H^{*}(Q \otimes \lambda)=0$ by Lemma \ref{repeat:demi1} applied  
to the $A_{2l-2}$ consisting of the first $2l-2$ simple roots.
For the cases where $j < 2l-2$, 
we can apply Lemma \ref{repeat:demi2}, 
also for the $A_{2l-2}$ consisting of the first $2l-2$ simple roots.
In that case, $a=j$, $b=2l-2$, $k=-j$
and so $k + b -a = 2l -2j - 2$, which is never $-1$.  
Also, clearly $-b-1 \leq k \leq -1$.
We therefore also have $H^{*}(Q \otimes \lambda)=0$.

Consequently, if we filter
$\wedge^j U_2$ by $B$-submodules
such that the quotients are one-dimensional, we deduce that 
$$H^{*} (S^{n-j} V_1^* \otimes \wedge^j U_2 \otimes \mu )=0$$
for  $1 \leq j \leq 2l-2$.  Hence Equation (\ref{step2}) follows.

\medskip

{\bf Step 3.}

In this step, we show that for all $i, n$
\begin{equation} \label{step3}
H^i (S^n V^* \otimes \mu) = 
H^{i} ( S^{n-l} V^* \otimes \mu + \omega_{2l} + \omega_{2l+1})
\end{equation}
for $\mu$ as in Step 2.

We take the Koszul resolution of the short exact sequence
$$0 \to U_1 \to V^{*} \to V_1^{*} \to 0$$ (this defines $U_1$)
and tensor it with $\mu$ arriving at 
\begin{equation} \label {cosi:1}
0 \to S^{n-2l+1} V^*  \otimes \wedge^{2l-1} U_1\otimes \mu \to \dots 
\to S^{n-j} V^* \otimes \wedge^j U_1 \otimes \mu \to \dots \to 
S^n V^* \otimes \mu \to S^n V_1^{*} \otimes \mu \to 0
\end{equation}

We first show that $H^*(S^{n-j} V^* \otimes \mu \otimes \lambda) = 0$
for any $\lambda$ appearing in $\wedge^j U_1$ for $j \neq 0, l$. 
The weights of $U_1$ are 
$$\al_k + \al_{k+1} + \dots + \al_{2l-1} + \al_{2l} + \al_{2l+1}$$ 
where $1 \leq k \leq 2l-1$.
%(and in particular, $\dim U_1 = 2l-1$).
If $\lambda$ is a weight of $\wedge^{j} U_1$, then as in the
previous steps we are quickly reduced to those $\lambda$ of the form
$$(1, 2, 3, \dots, j-2, j-1, j, \dots, j, j, j, j)$$
for $1 \leq j \leq 2l-1$.
Such a $\lambda$ satisfies 
$\langle \lambda, \al_{2l-1} \postcheck \rangle = -j$ 
with the exception of $j=2l-1$ where
$\langle \lambda, \al_{2l-2} \postcheck \rangle = -2l+2$.
The latter vanishing follows from Lemma \ref{repeat:demi1} applied to 
to the $A_{2l-1}$ consisting of the first $2l-1$ simple roots.
For the cases where $j < 2l-1$, 
we can apply Lemma \ref{repeat:demi2}, 
also for the $A_{2l-1}$ consisting of the first $2l-1$ simple roots.
In that case, $a=j$, $b=2l-1$, $k=-j$
and so $k + b -a = 2l -2j - 1$, which is $-1$
only when $j=l$.  Therefore, 
we deduce that 
$$H^{*} (S^{n-j} V^* \otimes \wedge^j U_1 \otimes \mu)=0$$
when $j \neq 0, l$.
And furthermore, 
$$H^{i} (S^{n-l} V^* \otimes \wedge^l U_1 \otimes \mu)= 
H^{i} (S^{n-l} V^* \otimes \lambda \otimes \mu),$$
where $\lambda = (1,2,3, \dots, l-1,l,\dots, l,l,l)$.
Now $S^{n-l} V^* \otimes \mu$ is stable under $P_{\al_m}$
for $1 \leq m \leq 2l-1$.  Hence $l-1$ applications of 
Proposition \ref{demazure} 
yields
$$H^{i} ( S^{n-l} V^* \otimes \lambda \otimes \mu)
=H^{i+l-1} (S^{n-l} V^* \otimes \mu + \omega_{2l} + \omega_{2l+1}).$$

By breaking Equation (\ref{cosi:1}) into short exact sequences 
and taking cohomology on $G/B$, we conclude that  
$$H^i (S^n V^* \otimes \mu) = 
H^{i} ( S^{n-l} V^* \otimes \mu + \omega_{2l} + \omega_{2l+1}),$$
where we are using 
$$H^{*} (S^{n} V_1^* \otimes \mu)=0$$
from Step 2.

\medskip

{\bf Step 4.}
We obtain the theorem by using Step 3 repeatedly, starting with $\mu = r \omega_m$
with $r$ in the prescribed range of the statement of the theorem.
After $-r$ steps we arrive at 
$$H^{i} (S^{n} V^* \otimes r \omega_{2l} ) = 
H^{i} (S^{n+rl} V^* \otimes -r \omega_{2l+1}),$$
for all $i, n$.
The proof is completed by using Step 1 and the symmetric version of Equation \ref{step1}
(obtained by applying an outer automorphism of $G$)
which gives 
$$H^{i} (S^{n+rl} V^* \otimes -r \omega_{2l+1}) =
H^i(S^{n+rl} \uni^*_{P'} \otimes -r \omega_{2l+1})$$
for all $i, n$.

\end{proof}

In what follows, we will use Theorem \ref{shift:A} in the more general
situation of Section 4 in \cite{sommers:normality}.  Similarly we 
can apply Theorem \ref{shift:D} in an analogous general situation.
Namely, suppose $G$ is of general type and 
$P$ is a parabolic subgroup of $G$ containing $B$ 
with Levi factor $L$ containing a simple factor of type $A_{2l}$.
Furthermore, suppose
this simple factor belongs to a Levi subgroup $L'$ of $G$ of type $D_{2l+1}$ and 
$[L, L'] \subset L'$.  Then the 
analog in $G$ of Theorem \ref{shift:D} holds 
just as the analog of Theorem \ref{shift:A} does 
in Proposition 6 in \cite{sommers:normality}.

\section{Main theorem}

For the rest of the paper $G$ is connected of type $D_{2l}$.  
We want to show that both nilpotent orbits in $\g$ with partition 
$(a^{2k}, b^2)$ for $a, b$ distinct even natural numbers with $a k + b = 2l$
(see \cite{col-mcg:nilp})
have normal closure.  Let $\orbit$ denote one of these two orbits.

Following the idea of \cite{sommers:normality}, we find a nilpotent orbit $\orbit'$ 
which we already know has normal closure and which contains $\orbit$ in its
closure.  
If we can show that the regular functions on $\orbit$ are naturally a quotient
of the regular functions on $\orbit'$, then it follows that $\orbit$ also has
normal closure. 
To that end we consider the nilpotent orbit $\orbit'$ in $\g$ with partition 
$\lambda = (a^{2k}, b+1, b-1)$.  

\begin{lem} 
The closure of $\orbit'$ is normal.
\end{lem}

\begin{proof}
The only minimal degenerations of $\orbit'$ in $\g$
are the two orbits with partition $\mu = (a^{2k}, b^2)$ (which together are one
orbit for the full orthogonal group of rank $2l$).
Hence by \cite{kraft-processi:normal2} the singularity of 
the closure of $\orbit'$ along the union of these two orbits 
is smoothly equivalent
to the singularity of the closure of the orbit with partition $(2)$ along 
the orbit with partition $(1,1)$ in type $A_1$
(we remove the first $2k$ rows from $\lambda$ and $\mu$, and then
remove the first $b-1$ columns from the resulting partitions).
Hence this is a singularity of type $A_1$ and so by 
\cite{kraft-processi:normal2},
$\orbit'$ has normal closure.
\end{proof}

\begin{lem}
The orbit $\orbit'$ is a Richardson orbit for any parabolic with Levi factor
of type 
$$\overbrace{A_{2k-1} \times \dots \times A_{2k-1}}^{\frac{a-b}{2} - 1}
 \times A_{2k} \times A_{2k} \times 
\overbrace{A_{2k+1} \times \dots \times A_{2k+1}}^{\frac{b}{2} - 1}.$$ 
%with birational moment map to $\g$. 
Any parabolic with Levi factor
of type 
$$\overbrace{A_{2k-1} \times \dots \times A_{2k-1}}^{\frac{a-b}{2}}
 \times 
\overbrace{A_{2k+1} \times \dots \times A_{2k+1}}^{\frac{b}{2}}$$
has Richardson orbit one or the other of the two nilpotent orbits
with partition $(a^{2k}, b^2)$.
\end{lem}

\begin{proof}
Both statements follow from Section 7 in \cite{col-mcg:nilp}. 
%The statement about the moment map follows since
%the equivariant fundamental group of $\orbit'$ is trivial. %cite
%The last statement follows from \cite{col-mcg:nilp} or
%just the fact that the parabolic listed here for $\orbit$ is its
%Jacobson-Morozov parabolic (corresponding to the non-negative
%eigenspaces of the semisimple part of $\mathfrak{sl}_2$-triple).
\end{proof}

%Let $P'$ denote the parabolic subgroup containing $B$ 
%in the previous lemma for $\orbit'$ 
% so that the Levi factors listed in the order of
%the lemma correspond to ordering on the simple roots.
%That is, the first Levi factor of type $A_{2k-1}$ 
%corresponds to the first $2k-1$ simple roots, etc.
%Let $\uni_{P'}$ be the Lie algebra of its unipotent radical. 

It will be convenient to represent parabolic subgroups
containing $B$ by the simple roots of $G$ which are 
{\bf not} simple roots of their Levi factors.  
Thus we can speak of such a parabolic subgroup as a subset of the numbers $1$ to $2l+1$,
with each number $i$ corresponding to the simple root $\al_i$.

%For example, let $\uni_{P'}$ be represented by
%$$\{ 2k, 4k, 6k, \dots, k(a-b-2), k(a-b)+1, k(a-b+2)+2, 
%k(a-b+4)+4, k(a-b+6)+6,  \dots,  ka +b = 2l \}.$$
%Then $\orbit'$ from the previous lemma is the Richardson orbit for $P'$
%(that is, $\orbit' \cap \uni_{P'}$ is dense in $\uni_{P'}$). 

%Let $P$ be the parabolic where $\uni_{P}$ is represented by
%$$\{ 2k, 4k, 6k, \dots, k(a-b-2), k(a-b), k(a-b+2)+2, 
%k(a-b+4)+4, k(a-b+6)+6,  \dots,  ka +b = 2l \}.$$
%Let $P_1$ be the parabolic 
%$$\{ 2k, 4k, 6k, \dots, k(a-b-2), k(a-b)+2, k(a-b+2)+2, 
%k(a-b+4)+4, k(a-b+6)+6,  \dots,  ka +b = 2l \}.$$

Set $d= a-b$ and let $P'$ be the parabolic represented by 
$$\{ 2k+1, 4k+2, 6k+2, \dots,  kd+2, k(d+2)+2, 
k(d+4)+4, k(d+6)+6,  \dots,  2l -2k-2,  2l \}$$
and let $P''$ be represented by
$$\{ 2k+1, 4k+2, 6k+2, \dots,  kd+2, k(d+2)+2, 
k(d+4)+4, k(d+6)+6,  \dots, 2l -2k-2,  2l-1 \},$$
so $P'$ are $P''$ are interchanged by an
outer automorphism of $D_{2l+1}$.  By the previous lemma 
$\orbit'$ is Richardson for both $P'$ and $P''$.
Let $P$ be the parabolic represented by 
$$\{ 2k, 4k+2, 6k+2, \dots,  kd+2, k(d+2)+2, 
k(d+4)+4, k(d+6)+6,  \dots, 2l -2k-2,  2l \}.$$
Then without loss of generality we can take $\orbit$ to be the Richardson orbit for $P$
(again by the previous lemma).

\begin{thm}
There is a short exact sequence
\begin{equation}
0 \to H^0( S^{n-2l-k(a-4)-1 }\uni^*_{P''} \otimes \nu) \to
H^0( S^{n} \uni^*_{P'}) \to
H^0( S^{n} \uni^*_{P}) \to
0,
\end{equation}
where $\nu = \omega_{4k+2}$ if $a>4$
%$$(1, 2, 3, \dots, 4k+1, \overbrace{4k+2, \dots, 4k+2}^{2l-4k-3}, 2k+1, 2k+1),$$
and $\nu = 2 \omega_{2l-1}$ if $a=4$ (and hence $b=2$).
%$$(1, 2, 3, \dots, 4k-2, 4k-1, 4k, 2k+1, 2k)$$
\end{thm}

\begin{proof}
We use two elements from the proof of Theorem \ref{shift:A} in \cite{sommers:functions}.
Let $P_1$ be the parabolic represented by 
$$\{ 2k+2, 4k+2, 6k+2, \dots,  kd+2, k(d+2)+2, 
k(d+4)+4, k(d+6)+6,  \dots, 2l -2k-2,  2l \}$$
and set $V = \uni_{P} \cap \uni_{P_1}$.
Then Step 1 of the proof of Theorem \ref{shift:A} 
(for a group of type $A_{4k+1}$ applied to the first
$4k+1$ simple roots of $G$) yields the isomorphism
$H^i( S^{n} \uni^*_{P}) = H^i( S^{n} V^*)$ for all $i,n$.
And Step 3 of the proof Theorem \ref{shift:A} 
yields the long exact sequence
$$ \dots \to H^i(S^{n-2k-1} \uni^*_{P'} \otimes \mu) \to 
H^i(S^{n} \uni^*_{P'}) \to 
H^i(S^{n} V^*) \to H^{i+1}(S^{n-2k-1} \uni^*_{P'} \otimes \mu) \to \dots$$
where $\mu$ equals
$$(1, 2, 3, \dots, 2k, 2k+1, 2k, \dots, 2, 1, 
\overbrace{0, 0, \dots, 0}^{2l-4k-1}).$$
This is obtained by taking the Koszul resolution of 
$$0 \to U \to \uni^*_{P'} \to V^* \to 0$$ (this defines $U$) and
simplifying the terms.

The remainder of the proof involves showing that
$$H^i(S^{n-2k-1} \uni^*_{P'} \otimes \mu) = 
H^i(S^{n-2l-k(a-4)-1} \uni^*_{P''} \otimes \nu)$$
for all $i,n$.

This is carried out by using Theorem \ref{shift:A} numerous times
(for $r=-1$ and the $l$ in that theorem equal to either $4k$ or $4k+1$ 
and $m'=2k$ or $2k+1$, respectively)
and Theorem \ref{shift:D} once (for $r=-2$ and 
the $l$ in that theorem equal to $k$).

After $\frac{a-b-2}{2}$ applications of 
Theorem \ref{shift:A} with $r=-1$, $l$ there equal to $4k$, and $m' = 2k$,
%applied to the leftmost $\frac{a-b-2}{2}$ factors, 
we have
$$H^i(S^{n-2k-1} \uni^*_{P'} \otimes \mu) =
H^i(S^{n -k(a-b) -1} Q^*_1 \otimes \mu_1)$$
where $\mu_1$ equals
$$(1, 2, 3, \dots, 2k, \overbrace{2k+1, \dots, 2k+1}^{k(a-b-2)+1}, 2k, \dots, 2, 1, 
\overbrace{0, 0, \dots, 0}^{k(b -2)+b -1}),$$
and $Q_1$ is the Lie algebra of the unipotent radical of  
$$\{ 2k+1, 4k+1, 6k+1, \dots, k(d-2)+1, kd+1, k(d+2)+2, 
k(d+4)+4, k(d+6)+6,  \dots, 2l-2k-2, 2l \}.$$

Next, 
we apply Theorem \ref{shift:A}
$\frac{b-2}{2}$ more times with $r=-1$, $l$ there equal to $4k+1$, and $m' = 2k+1$,
%applied to the rightmost $\frac{b-2}{2}$ factors, 
to obtain 
$$H^i(S^{n -k(a-b) -1} Q^*_1 \otimes \mu_1) = 
H^i(S^{n -ka +2k -b/2} Q^*_2 \otimes \mu_2)$$
where $\mu_2$ equals
$$(1, 2, 3, \dots, 2k, \overbrace{2k+1, \dots, 2k+1}^{2l-4k-1}, 
2k, \dots, 2, 1, 0),$$
and $Q_2$ is the Lie algebra of the unipotent radical of  
$$\{ 2k+1, 4k+1, 6k+1, \dots, k(d-2)+1, kd+1, k(d+2)+3, 
k(d+4)+5, k(d+6)+7,  \dots, 2l -2k-1,  2l \}.$$

Next, we use Theorem \ref{shift:D}
with $r=-2$ for the case $D_{2k+1}$ applied to the simple
roots $\al_i$ of $G$ with $2l-2k \leq i \leq 2l$. 
This yields 
$$H^i(S^{n -ka +2k -b/2} Q^*_2 \otimes \mu_2) = 
H^i(S^{n -ka -b/2} Q^*_3 \otimes \mu_3)$$
where $\mu_3$ equals
$$(1, 2, 3, \dots, 2k, \overbrace{2k+1, \dots, 2k+1}^{2l-4k-1}, 
2k+2, 2k+3, 2k+4, \dots, 4k, 2k+1, 2k),$$
and $Q_3$ is the Lie algebra of the unipotent radical of  
$$\{ 2k+1, 4k+1, 6k+1, \dots, k(d-2)+1, kd + 1, k(d+2)+3, 
k(d+4)+5, k(d+6)+7,  \dots, 2l -2k-1,  2l-1 \}.$$

If $2l-4k-1=1$, which is the case if and only if $a=4$ and $b=2$, 
we have $\mu_3 = 2 \omega_{2l-1}$ and the latter parabolic subgroup is $P''$.

On the other hand, if $a>4$,
we continue by using Theorem \ref{shift:A}
another $\frac{b-2}{2}$ times followed by another $\frac{a-b-2}{2}$
times (in reverse of how we have just used it).  The result is that 
$$H^i(S^{n-ka-b/2-2 } Q^*_3 \otimes \mu_3)=
H^i(S^{n-2ka+4k-b-1 } Q^*_4 \otimes \mu_4)$$
where $\mu_4$ equals
$$(1, 2, 3, \dots, 4k+1, \overbrace{4k+2, \dots, 4k+2}^{2l-4k-3}, 2k+1, 2k+1),$$
and $Q_4$ is the Lie algebra of the unipotent radical of  
$$\{ 2k+1, 4k+2, 6k+2, \dots, k(d-2)+2, kd+2, k(d+2)+2, 
k(d+4)+4, k(d+6)+6,  \dots, 2l -2k-2,  2l-1 \}.$$
The latter parabolic is exactly $P''$ and $\mu_4 = \omega_{4k+2}$. 
Furthermore, $n\!-\!2ka + 4k\!-\!b\!-\!1 = n\!-\!2l\!-\!ka + 4k\!-\!1$ since $ak+b = 2l$.

Hence when $a=4$ or $a>4$, we have shown that 
$$H^i(S^{n-2k-1} \uni^*_{P'} \otimes \mu) = 
H^i(S^{n-2l-k(a-4)-1} \uni^*_{P''} \otimes \nu)$$
for all $i,n.$  
We finish the proof by observing that $\nu$ extends to a character of $P''$ and it is dominant.   
Hence $H^i(S^{n-2l-k(a-4)-1} \uni^*_{P''} \otimes \nu) = 0$
for $i>0$ as in \cite{broer:vanishing}.
Similarly,
$H^i( S^{n} \uni^*_{P'})=0$ and $H^i( S^{n} \uni^*_{P})=0$ for $i>0$ and the proof is complete. 
\end{proof}  

\begin{cor}
The closure of $\orbit$ is normal.
\end{cor}

\begin{proof}
We only need to note that the functions 
of degree $n$ on $\orbit'$ (and also its closure since
the closure is normal) as a $G$-module are isomorphic to $H^0( S^{n} \uni^*_{P'})$.
This follows since $\orbit'$ has trivial $G$-equivariant fundamental group when $G$ is adjoint
(see \cite{col-mcg:nilp}).  Hence the moment map determined by $P'$ must be birational.
Thus the short exact sequence of the theorem together with the discussion 
in Section 3 of \cite{sommers:normality} yields the result.
\end{proof}

%\begin{rmk}
%Other very evens have an analog of the theorem, but in those
%cases the analogous orbit $\orbit'$ does  not have normal closure (that is,
%Lemma  does not hold).  This explains why the corollary does not go
%through in those cases and indeed from \cite{kraft-processi:normal2} we know that
%$\orbit$ does not have normal closure.
%\end{rmk}

\bibliography{sommers_normal_D}
\bibliographystyle{pnaplain}
\end{document}